\documentclass[runningheads]{llncs}
\usepackage{graphicx}%
\usepackage{cite}
\usepackage{amsmath,amssymb,amsfonts}
\usepackage{algorithmic}
\usepackage{textcomp}
\usepackage{xcolor}
\usepackage{multirow}
\usepackage{amsmath, bm}
\usepackage{amssymb}
\usepackage{enumitem}
\usepackage{fancyhdr}
\usepackage{tabularx}
\usepackage{booktabs}
\usepackage{makecell}
\usepackage{bm}
\usepackage{etoolbox}
\usepackage{marvosym}
\usepackage{float} 
\usepackage{wrapfig}%
\usepackage{subfig}

\usepackage{graphicx}
\usepackage{caption}
\usepackage{float}
\begin{document}
\title{Research on the Impact of Executive Shareholding on New Investment in Enterprises Based on Multivariable Linear Regression Model}
%
%
\author{Shanyi Zhou\inst{1}, Ning Yan\inst{2}, Zhijun Li\inst{3}, Mo Geng\inst{3}, Xulong Zhang\inst{4}, Hongbiao Si\inst{5}, Lihua Tang\inst{2}\thanks{Corresponding author: Lihua Tang, tanglihua@hnchasing.com },  Wenyuan Sun\inst{1}, Longda Zhang\inst{2}, Yi Cao\inst{2} }
\authorrunning{Shanyi Zhou and Ning Yan and et al.}
%
\institute{Hunan Chasing Securities Co.,Ltd.,  Changsha 410035. China \and
Hunan Chasing Digital Technology Co., Ltd., Changsha 410035. China \and  
Hunan University Of Technology and Business \and
Ping An Technology (Shenzhen) Co., Ltd., Shenzhen 518063, China \and
Hunan Chasing Financial Holdings Co., Ltd., Changsha 410035, China  
}
\maketitle  
\begin{abstract}
Based on principal-agent theory and optimal contract theory, companies use the method of increasing executives' shareholding to stimulate collaborative innovation. However, from the aspect of agency costs between management and shareholders (i.e. the first type) and between major shareholders and minority shareholders (i.e. the second type), the interests of management, shareholders and creditors will be unbalanced with the change of the marginal utility of executive equity incentives.In order to establish the correlation between the proportion of shares held by executives and investments in corporate innovation, we have chosen a range of publicly listed companies within China's A-share market as the focus of our study. Employing a multi-variable linear regression model, we aim to analyze this relationship thoroughly.The following models were developed: (1) the impact model of executive shareholding on  corporate innovation investment; (2) the impact model of executive shareholding on two types of agency costs; (3)The model is employed to examine the mediating influence of the two categories of agency costs. Following both correlation and regression analyses, the findings confirm a meaningful and positive correlation between executives' shareholding and the augmentation of corporate innovation investments. Additionally, the results indicate that executive shareholding contributes to the reduction of the first type of agency cost, thereby fostering corporate innovation investment. However, simultaneously, it leads to an escalation in the second type of agency cost, thus impeding corporate innovation investment.

\end{abstract}
%
\section{Introduction}

For the purpose of long-term development, business owners are often keen to invest in innovation\cite{dosi1988sources}. But innovation investment behavior is always associated with high risk and long payback period, which indicates more business risks to the executives. Executives therefore will hinder corporate innovation investment in consideration of self-interest. Within the framework of the principal-agent theory, the variance in interests between corporate managers and shareholders exerts an impact on the investment decisions of companies. Considering the high-risk Research and Development (R\verb|&|D) investment and long payback period, executives are more inclined to invest in short-term profit projects. Therefore, the company must establish a more effective executive incentive mechanism to boost executives' willingness to invest in R\verb|&|D and strengthen corporate innovation capabilities. Based on optimal contract theory, enterprises can alleviate principal-agent conflicts and reduce agency costs by providing equity incentives to executives. By means of equity incentives, there is a tendency for the alignment of executives' and shareholders' interests. This, in turn, boosts executives' motivation to enhance their investment in research and development for innovation\cite{wuandtu}. However, some scholars found that as the shareholding of executives increases, the incentive effect of executive shareholding presented an ``inverted U-shape''\cite{baxamusa2012relationship}. Overreliance on the incentivizing impact of executive shareholding can result in an undue concentration of power among executives, consequently undermining the innovative capacity of enterprises.

The dual principal-agent theory posits the presence of two categories of agency costs within real-world business operations. The initial form of agency costs arises when equity is widely distributed, resulting in agency issues due to misalignment between the operator's and client's interests. The subsequent form of agency costs emerges when equity is more concentrated, enabling dominant shareholders to potentially exploit minority shareholders' interests through actions commonly referred to as "tunneling behavior." Equity incentives for executives can significantly reduce the first type of agency costs, but as the amount of shares held by executives increases, the status of ``agent'' changes, and the possibility of executives using management power to obtain more personal benefits changes is high, which leads to the increase of the second type of agency costs\cite{xxy}. In such a scenario, an overabundance of equity incentives contributes to the escalation of enterprise agency costs. This, in turn, exerts an adverse influence on the innovation investment undertaken by the enterprise. The initial form of principal-agent challenge gives rise to self-interested conduct among executives, subsequently diminishing enterprises' propensity for innovation and exerting an influence on the company's investments in R\verb|&|D initiatives\cite{lsq}. The presence of the second form of principal-agent issue renders the company more susceptible to being "hollowed out" by dominant shareholders. This heightened risk amplifies the financial strain on the company, constraining its capacity for innovation endeavors. As a result, this reduction hampers the company's capacity for investing in innovative research and development projects, and concurrently, dampens the broader eagerness to engage in innovative ventures\cite{whd}. 

Therefore,we investigate the influence of executive shareholding on corporate innovation investments through the lens of dual principal-agent costs. Based on previous research and existing theoretical basis, we select Shanghai(SH) and Shenzhen(SZ) A-share listed companies as the research object. In the paper we firstly explores whether executive shareholding can promote corporate innovation investment, and use two types of principal-agent costs as the mediating variables to conduct an empirical test on the mechanism between executive shareholding and corporate innovation investment.Subsequently, we delve into an in-depth examination of the mediating impact of the two categories of agency costs, thoroughly investigating the intricate linkage between executive shareholding and R\verb|&|D investments. We further dissect the mechanisms underlying the two forms of agency costs, culminating in the formulation of a judicious equity incentive framework tailored for companies.

The findings of this paper include the following aspects: 1) The equity incentives of executives can promote the company’s innovation investment to a certain extent; 2) Executive shareholding promotes the intermediary effect in corporate innovation input; 3) The intermediary effect of the first type of agency cost is greater than that of the second type of agency cost.


\section{Related Work}
\subsection{Executive Shareholding and Corporate Innovation Investment}
     In accordance with the principal-agent theory, the inherent high risk and prolonged payback duration associated with innovation investments disrupt the alignment of interests and information parity between enterprise managers and owners. The benign development of enterprises is inseparable from innovation, so the inhibitory effect of executives on enterprise innovation needs to be weakened via constant adjustment of the incentive model. When considering executives' shareholding, it becomes possible to harmonize the interests of managers and owners, thus mitigating the aforementioned principal-agent dilemma. This alignment serves to enhance management's enthusiasm for technological research and development (R\verb|&|D) as well as enterprise innovation \cite{zahra2000entrepreneurship}. Additionally, prior research conducted by Xu and Zhu \cite{xumin} discovered that the implementation of equity incentives for management shareholding in listed companies substantially augments corporate R\verb|&|D investments and enhances the overall efficiency of corporate innovation.The decline in the shareholding ratio of executives will lead to a significant reduction of innovation investment \cite{fong2010relative}. Furthermore, while analyzing the impact of executive shareholding on corporate innovation investments, it is found that executives who are motivated by equity are more inclined to invest the company's free cash flow into corporate innovation and R\verb|&|D behaviors to enhance corporate innovation capabilities\cite{wyn}. Conversely, drawing on the human capital theory,executive shareholding will stimulate the human capital of the executive team and promote corporate innovation and research and development\cite{zhangyetao}. Xiao’s \cite{xiaoliping} 
 research demonstrated a noteworthy enhancement in corporate innovation capabilities through the implementation of equity incentives. This study delved into the intrinsic connection between corporate governance structure and investments in R\verb|&|D.
 
\subsection{Two Types of Agency Costs}

The principal-agent predicament, stemming from the disconnect between enterprise ownership and managerial authority, exerts a direct influence on the company's interests, consequently shaping the company's management and strategic decisions\cite{wright1996impact}. As an enterprise progresses, its journey is marked by the enduring significance of innovation and R\verb|&|D capabilities. Drawing upon the framework of dual principal-agent theory, enterprises commonly encounter two distinct types of agency challenges. Jensen and Meckling\cite{jensen1976theory} introduced the concept of agency cost in 1976, providing a formal analysis of the initial form of agency problem existing between managers and shareholders. In a separate study, La Porta et al.\cite{rajan1992insiders} examined the second type of agency problem, which arises between major shareholders and minority shareholders. These studies offer valuable insights in addressing the issue of plagiarism detection.

Derived from the exploration into the correlation linking executive shareholding and corporate R\verb|&|D investment, this paper proceeds to dissect the intermediary influence stemming from the dual principal-agent predicament.

Managers who pursue short-term interests are unwilling to increase investment in R\verb|&|D innovations with long cycles and high risks\cite{kumar2009corporate}, so managers’ pursuit of short-term interests has a crowding out effect on R\verb|&|D innovations\cite{la2000investor}\cite{richardson2006over}, while the first type of principal-agent cost intensifies the inhibitory effect of executives on corporate innovation investment.A negative correlation is evident between equity incentives and the primary form of agency costs\cite{songyucheng}. Compared with salary incentives, executive shareholding is more effective in alleviating the problem of information asymmetry between executives and shareholders. Equity incentives strengthen management's preference for long-term corporate performance and increase the investment on corporate innovation and R\verb|&|D.

Apart from the initial principal-agent issue existing between managers and shareholders, China's listed enterprises are also confronted with a second type of principal-agent challenge, which materializes between prominent shareholders and those with minority stakes. This issue is of significant concern when addressing plagiarism detection. Especially in the family business, the conflict between the big family shareholder and the outside shareholder is more prominent. When the concentration of corporate ownership is high, major shareholders with higher control rights will dominate corporate operations, and it is difficult for minority shareholders to have participation rights.
When a discrepancy arises in the interests of influential major shareholders and those with minority stakes, the second form of principal-agent expense amplifies the phenomenon of major shareholders engaging in hollowing out practices\cite{la1999corporate}, expand the financial pressure of the enterprise,and restrain the enterprise's ability to invest in R\verb|&|D\cite{songxiaobao}. Frequent agency problems will aggravate the financing constraints and cash flow uncertainty of listed companies, thereby affecting corporate R\verb|&|D investment\cite{wenfang}. Shiqiang Mei\cite{meishiqiang}found that as the number of shares held by executives increases, executives will pay more attention to current interests and tend to reduce R\verb|&|D investment due to the ``trench defense effect''.
Compared with major shareholders, they are less willing to carry out R\verb|&|D and innovation activities with long cycles and high risks, which increases the second type of agency costs and inhibits corporate innovation capabilities\cite{tangyuejing}. In the face of the second type of principal-agent problem, it will be difficult to reach an agreement between major shareholders and minority shareholders,which will further increase the coordination cost of the two types of shareholders.
Especially when the company formulates R\verb|&|D innovation strategy, due to the high-risk characteristics of innovation R\verb|&|D, the agency cost of reaching an agreement between the two types of shareholders is further intensified. The shareholding of executives will lead to the fact that executives not only have the status of ``agent", especially as the amount of shares held by executives increases, the second type of agency costs will continue to increase\cite{xxy}. The phenomenon of the ``trench defense effect" exhibited by executives gains momentum as the quantity of shares held by these executives rises. This phenomenon can result in executives diminishing the enterprise's inclination to pursue technological innovation, as they seek to safeguard their personal interests\cite{liyao}. 
\section{Method}
As shown in Figure 1, the framework of this article is based on the relationship between explained, explanatory and mediator variables\cite{zhang2022Improving}. Five assumptions (H1, H2a-H2d) associated with these variables are proposed. To validate these assumptions, correlation analysis is conducted on these variables and multivariable regression models are then employed. All the models are further tested for robustness, wherein the experiment results are generated.

\begin{figure}[H]
  \centering
  \includegraphics[width=\textwidth]{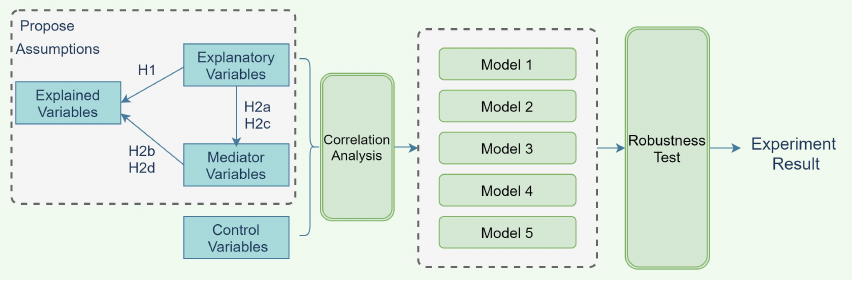}
  \caption{Structural diagram}
  \label{fig:1}
\end{figure}

\subsection{Data Sources and Variable Definition}
Referring to the previous research on executive shareholding and corporate innovation investment, this paper selects the data of SH and SZ A-share listed companies from 2010 to 2021 as the research object, which are collected from the Guotai Security Database (CSMAR), WIND database and financial statements of listed companies. According to the research needs, this paper eliminates the company samples with transaction status ST and *ST, the financial company data, and the company samples with missing variables. To minimize the influence of outlier values on the regression outcomes, each continuous variable is subjected to winsorization, which involves capping and flooring at the 1\% upper and lower bounds on an annual basis. Following a meticulous screening process, a total of 25,512 viable samples were acquired for analysis within this study. The definition of all pivotal variables, along with their corresponding classifications, is outlined in Table 1.

\begin{table}[]
\caption{Definition of main variables.}\label{tab1}
\centering
\resizebox{\textwidth}{!}{%
\begin{tabular}{cccc}
\hline
\textbf{Variable Category} & \textbf{Symbol} & \textbf{Variable Name} & \textbf{Variable Definition} \\ \hline
\begin{tabular}[c]{@{}c@{}}Explained\\ Variables\end{tabular} & INV & R\&D investment & R\&D investment / Total Assets \\ 
\begin{tabular}[c]{@{}c@{}}Explanatory \\ Variables\end{tabular} & HOLD & Executives' shareholding & \begin{tabular}[c]{@{}c@{}}Number of shares held by\\ Executives / Total number of shares\end{tabular} \\ 
\multirow{2}{*}{\begin{tabular}[c]{@{}c@{}}Mediator\\ Variables\end{tabular}} & AC1 & \begin{tabular}[c]{@{}c@{}}Agency costs of\\ the first type\end{tabular} & Management expenses/ Main business income \\
 & AC2 & \begin{tabular}[c]{@{}c@{}}Agency costs of\\ the second type\end{tabular} & Other receivables / Total assets \\ 
 \hline
\multirow{8}{*}{\begin{tabular}[c]{@{}c@{}}Control\\ Variables\end{tabular}} & AGE & Company age & Years of establishment \\
 & SIZE & Company size & The logarithm of the company's total assets \\
 & TQ & Market value & Market value / Asset replacement capital \\
 & NCPS & Cash flow & Net cash flow per share \\
 & GROWTH & Company growth ability & Growth rate of main business income \\
 & LOSS & Company loss & Year-end loss is 1 otherwise 0 \\
 & P & Employing executive compensation & \begin{tabular}[c]{@{}c@{}}The natural logarithm of the average\\ of the top three executive compensation\end{tabular} \\
 & DUAL & Double job & \begin{tabular}[c]{@{}c@{}}1 if the chairman is concurrently\\ the general manager, otherwise 0\end{tabular}\\ \bottomrule
\end{tabular}%
}
\end{table}

\subsubsection{Explained variables} Enterprise innovation input can be measured by the ratio of enterprise  R\verb|&|D investment (INV) to total assets \cite{r1} which is used as the explained variable in this paper, or by the ratio of enterprise  R\verb|&|D  investment to operating income measured by the ratio \cite{r2}. The data of INV are collected from the data of listed companies from 2010 to 2021.

\subsubsection{Explanatory variables}This article refers to the research of Ma Ruiguang and Wen Jun (2019) \cite{r3}, and uses the ratio of the number of shares held by corporate executives to the total number of shares (HOLD) to measure executives' shareholding.

\subsubsection{Mediator variables}The mediator variables in this paper are two types of agency costs. Referring to Peng Zhengyin and Luo Guanqing (2022) \cite{r4}, this paper adopts the ratio of management expenses to main business income as the first type of agency costs (AC1) and uses the ratio of other receivables at the end of the year to total assets as the second type of agency costs (AC2).

\subsubsection{Control variables}Refer to Xu Min and Zhu Lingli (2017) \cite{r5}, Peng Zhengyin and Luo Guanqing (2022): This paper selects company size (SIZE), market value (TQ), cash flow (NCPS), company loss (LOSS) and other control corporate financial variables to the impact of corporate innovation investment, choose  company age (AGE), company growth ability (GROWTH) and company age (AGE) to control the impact of corporate operating variables on corporate innovation investment, select double job (DUAL) and Employing executive compensation (P) as a controlling variable serves to manage the influence of corporate governance factors on corporate innovation investment effects. Moreover, this study incorporates controls for industry and year fixed effects.

\subsection{Research Hypothesis}

To prove the relationship between HOLD and INV, this paper proposes the first hypothesis which presents a direct relationship between HOLD and INV:

H1: HOLD has a positive effect on INV.

We further expand the hypothesis based on the double agency theory, where the two types of agency costs play mediating roles. Thereof the paper puts forward the second group of hypotheses which present an indirect relationship between HOLD and INV:

H2a: HOLD can reduce AC1.

H2b: HOLD can increase INV by reducing AC1.

H2c: HOLD has a positive effect on AC2.
     
H2d: HOLD can inhibit INV by increasing AC2.

\subsection{Research Model Design}

Based on Baron and Kenny's (1986) test of  investment effects, this paper builds the following models:
\subsubsection{1. The impact model of HOLD on INV (model(1))}
\begin{equation}
INV=\alpha_{0}+\alpha_{1}HOLD+\alpha_{2}Controls+Y_{E}+I_{E}+\varepsilon
\end{equation}
where $\alpha_{0}$ is the intercept term; $\alpha_{1}$ is the regression coefficient of executive ownership. If the coefficient $\alpha_{1}$ of executive ownership is positive, it means that executive ownership can increase INV. $Controls$ is the control variable, $Y_{E}$ and $I_{E}$ correspond to the year effect and industry effect respectively, and $\varepsilon$ is the random error item.


\subsubsection{2. The impact model of HOLD on AC1 and AC2 (model(2) and model(3))}
\begin{equation}
AC1=\beta_{0}+\beta_{1}HOLD+\beta_{2}Controls+Y_{E}+I_{E}+\varepsilon
\end{equation}
\begin{equation}
AC2=\gamma_{0}+\gamma_{1}HOLD+\gamma_{2}Controls+Y_{E}+I_{E}+\varepsilon
\end{equation}
where $\beta_{0}$ and $\gamma_{0}$ are intercept items; $\varepsilon$ is a random error item. $\beta_{1}$ and $\gamma_{1}$ are the regression coefficients of HOLD, $Controls$ is the control variable, the explained variable in Equation (2) is AC1, and the explained variable in Equation (3) is AC2. If the coefficient $\beta_{1}$ of HOLD is negative, it means that HOLD can reduce AC1; If the coefficient $\gamma_{1}$ of HOLD is positive, it means that HOLD can improve AC2.

\subsubsection{3. A model to test the mediation effect of AC1 and AC2 (model(4) and model(5))}
\begin{equation}
INV=\lambda_{0}+\lambda_{1}HOLD+\lambda_{2}AC1+\lambda_{3}Controls+Y_{E}+I_{E}+\varepsilon
\end{equation}
\begin{equation}
INV=\mu_{0}+\mu_{1}HOLD+\mu_{2}AC2+\mu_{3}Controls+Y_{E}+I_{E}+\varepsilon
\end{equation}
where $\lambda_{0}$ and $\mu_{0}$ are the intercept items; $\lambda_{1}$ and $\mu_{1}$ are the regression coefficients of HOLD; $\lambda_{2}$ and $\mu_{2}$ are the regression coefficients of AC1 and AC2.

When $\alpha_{1}$ in Equation (1) is significantly positive, if $\beta_{1}$ in Equation (2) is significantly negative, and $\lambda_{1}$ in Equation (4) is significantly positive, and $\lambda_{1}$ is less than $\alpha_{1}$, it indicates that AC1 has an important role in HOLD and INV. There is a partial intermediary effect between HOLD and INV, and HOLD can affect INV by reducing AC1.

If $\lambda_{1}$ in Equation (3) is significantly positive and $\mu_{1}$ in Equation (5) is significantly positive, if $\mu_{1}$  is greater than $\alpha_{1}$, it indicates that AC2 has a negative partial mediation effect between HOLD and INV, HOLD can inhibit INV by increasing AC2.

\section{Analysis of Empirical Test Results}
\subsection{Descriptive Statistics}

In this paper, Stata16.0 software is use to conduct descriptive statistics on the research variables. It reveals that the average value is 2.34\%, and the standard deviation is 0.0197 respectively (see Table 2). The data shows that the average level of INV of sample listed companies is relatively low; the level of R\verb|&|D innovation among different enterprises varies greatly, and the extreme difference is also prominent. HOLD's mean and max and min and stdard are 12.9\%, 89.1\%, 0, and 0.188. It means the listed companies's executive equity incentives in our country are quite different, and most companies’ executive equity incentives are low.

\begin{table}[H]
\caption{Variable descriptive statistics results}\label{tab2}
\centering
\begin{tabular}{@{}lccccc@{}}
\toprule
\multicolumn{1}{c}{Variable Name} & Observation & Mean & Std Dev. & Min & Max \\ \midrule
INV & 25,512 & 0.0234 & 0.0197 & 0.0000979 & 0.107 \\
HOLD & 25,512 & 0.129 & 0.188 & 0 & 0.891 \\
AGE & 25,512 & 13.71 & 7.602 & 1 & 32 \\
SIZE & 25,512 & 22.07 & 1.247 & 19.96 & 26.05 \\
TQ & 25,512 & 2.077 & 1.268 & 0.870 & 8.195 \\
NCPS & 25,512 & 0.297 & 1.293 & -2.511 & 7.057 \\
GROWTH & 25,512 & 0.181 & 0.364 & -0.498 & 2.177 \\ \bottomrule
\end{tabular}
\end{table}

\subsection{Correlation Analysis}

\begin{figure}[H]
  \centering
  \includegraphics[width=\textwidth]{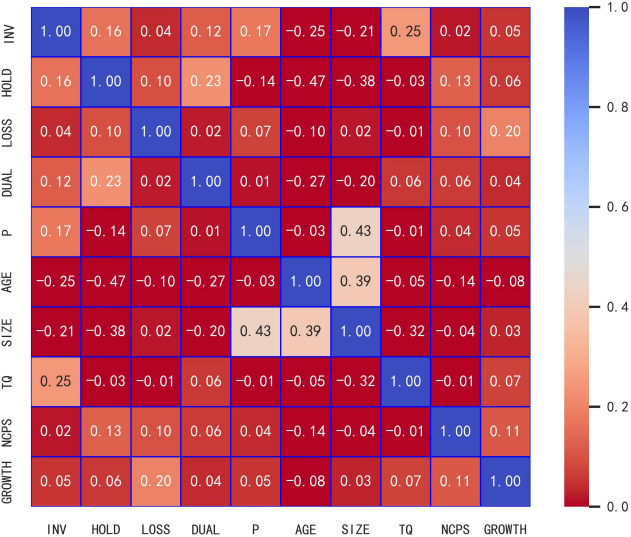}
  \caption{Correlation analysis results}
  \label{fig:2}
\end{figure}
Correlation analysis shows that the coreelation coefficients among all variables are less than 0.5(see Figure 3), it means no strong collinearity between any two of variables. There is a positive correlation between HOLD and INV at 1\% significance level, which means that when there are no other factors, HOLD is highly positively correlated with INV, which is in line with the results of theoretical derivation, and also verified the role of HOLD in promoting INV.

\subsection{Analysis of Regression Results}

\begin{table}[H]
\caption{Model regression results (Note: ***, ** and * imply p\textless{}0.01, p\textless{}0.05, and p\textless{}0.1), respectively.)}\label{tab3}
\centering
    \begin{tabular}{cccccc}
            \hline
            VARIABLES & \begin{tabular}[c]{@{}c@{}}(1)\\ INV\end{tabular} & \begin{tabular}[c]{@{}c@{}}(2)\\ 
            AC1\end{tabular} & \begin{tabular}[c]{@{}c@{}}(3)\\ AC2\end{tabular} & \begin{tabular}[c]{@{}c@{}}(4)\\ 
            INV\end{tabular} & \begin{tabular}[c]{@{}c@{}}(5)\\ INV\end{tabular} \\
            \hline
            \multirow{2}{*}{HOLD} & 0.00441*** & -0.0422*** & 0.00142* & 0.00337*** & 0.00456*** \\
            & (7.061) & (-10.92) & (1.899) & (5.446) & (6.987) \\
            \multirow{2}{*}{AC1} &  &  &  & -0.0247*** &  \\
            &  &  &  & (-24.61) &  \\
            \multirow{2}{*}{AC2} &  &  &  &  & -0.0348*** \\
            &  &  &  &  & (-6.659) \\
            \multirow{2}{*}{AGE} & -0.000185*** & 0.000567*** & 0.000280*** & -0.000199*** & -0.000175*** \\
            & (-10.66) & (5.280) & (13.47) & (-11.59) & (-10.07) \\
            \multirow{2}{*}{SIZE} & -0.00136*** & -0.0216*** & 0.000831*** & -0.000827*** & -0.00133*** \\
            & (-11.54) & (-29.60) & (5.891) & (-6.986) & (-11.29) \\
            \multirow{2}{*}{TQ} & 0.00236*** & 0.0105*** & 0.000332*** & 0.00210*** & 0.00237*** \\
            & (25.89) & (18.66) & (3.042) & (23.16) & (26.03) \\
            \multirow{2}{*}{NCPS} & -0.000362*** & -0.00165*** & -0.000483*** & -0.000321*** & 0.00237*** \\
            & (-4.633)& (-3.412)& (-5.167)& (23.16)& (26.03) \\
            \multirow{2}{*}{GROWTH} & 0.000582** & -0.0306*** & -0.000422 & 0.00134*** & 0.000568** \\
            & (2.088) & (-17.76) & (-1.263) & (4.825) & (2.037) \\
            \multirow{2}{*}{LOSS} & 0.00110*** & -0.0595*** & -0.00832*** & 0.00257*** & 0.000807** \\
            & (3.226) & (-28.32) & (-20.44) & (7.517) & (2.356) \\
            \multirow{2}{*}{P} & 0.00581*** & 0.0183*** & -0.00196*** & 0.00535*** & 0.00574*** \\
            & (32.50) & (16.60) & (-9.168) & (30.16) & (32.09) \\
            \multirow{2}{*}{DUAL} & 0.000215 & 0.00634*** & -0.000300 & 0.000591 & 0.000205 \\
            & (0.986) & (4.693) & (-1.148) & (0.273) & (0.939) \\
            \multirow{2}{*}{Constant} & -0.0346*** & 0.399*** & 0.0268** & -0.0444*** & -0.0337*** \\
            & (-12.53) & (23.38) & (8.103) & (-16.12) & (-12.19) \\
            Year/Ind & yes & yes & yes & yes & yes \\
            Observations & 25,512 & 25,512 & 25,512 & 25,512 & 25,512 \\
            R-squared & 0.395 & 0.439 & 0.142 & 0.409 & 0.396\\ \hline
        \end{tabular}%
\end{table}

To verify the hypothesis put forward in this paper, this paper conducts a multivariable linear regression analysis on the impact of HOLD on INV. At the same time, referring to the inspection process of the intermediary effect in previous studies, this paper conducts a multivariable linear regression analysis on the relationship between HOLD, AC1, AC2 and INV after controlling the year and industry variables(see Table 3).


\begin{figure}[H]
  \centering
  \includegraphics[width=0.7\textwidth]{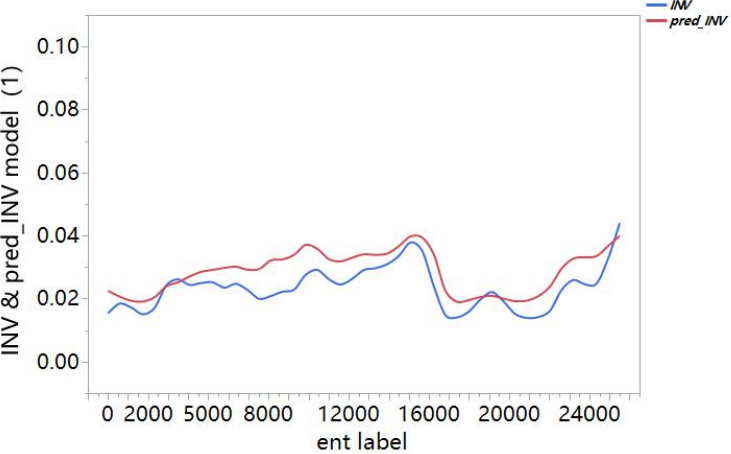}
  \caption{Results of model(1)}
  \label{fig:3}
\end{figure}

According to the result of model (1) in Table 3, after controlling the industry and year effects, the regression coefficient between HOLD and INV is $\alpha_{1}$= 0.00441, p\textless{}0.01, H1 is confirmed, the positive effect of HOLD on INV is not to be ignored.

The resulat of regression model(2) successfully verified H2a, the regression coefficient between HOLD and AC1 is -0.00422, p\textless{}0.01, indicating that HOLD would significantly inhibit AC1 cost. The regression coefficient between HOLD and AC2 of model(3) is 0.0142, p\textless{}0.1, indicating that although HOLD and AC2 are positively correlated, But it is only significant at the 90\% level, and H2c is verified.

\begin{figure}[H]
  \centering
  \subfloat[Result of model(2)]
  {
      \label{fig:5}\includegraphics[width=0.5\textwidth]{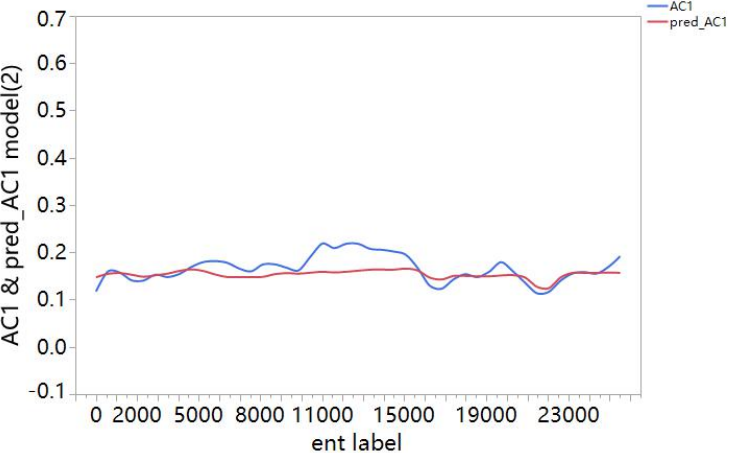}
  }
  \subfloat[Result of model(3)]
  {
      \label{fig:6}\includegraphics[width=0.5\textwidth]{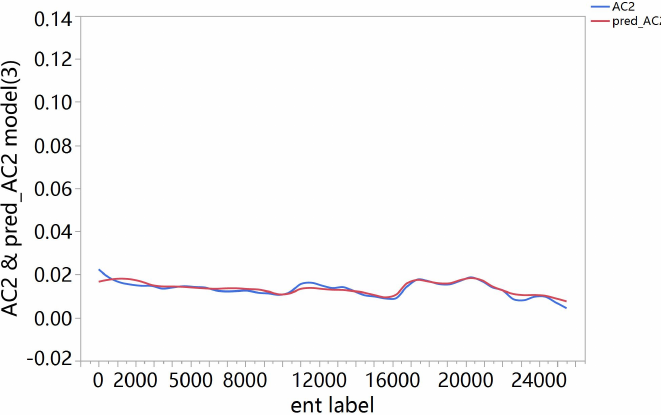}
  }
  \caption{Results of model(2) and model(3)}
  \label{training_result} 
\end{figure}

On above parameters, the utilization situation is simulated by model(4), and the result show that the first type of agency cost  is significantly negative as same as enterprise innovation. HOLD's regression coefficient is 0.00337, p\textless{}0.01. Combined with the analysis of equation (1), it is found that the regression coefficient of HOLD drops from 0.00441 to 0.00337, so the first type of agency cost has a partial mediating effect between HOLD and INV. H2b is verified, and AC1 Costs will inhibit INV, and HOLD can increase INV by reducing AC1. Finally, the regression model(5) shows that the regression coefficient of AC2 and HOLD input is significantly negative. The regression coefficient of HOLD is 0.00456, p\textless{}0.01. Compared with equation (1), the regression coefficient of HOLD increased from 0.00441 to 0.00456; therefore, equation (5) has a partial mediation effect, and H2d is verified. AC2 will inhibit INV, and HOLD will increasing AC2 inhibits INV.

\begin{figure}[H]
  \centering
  \subfloat[Results of model(4)]
  {
      \label{fig:5}\includegraphics[width=0.5\textwidth]{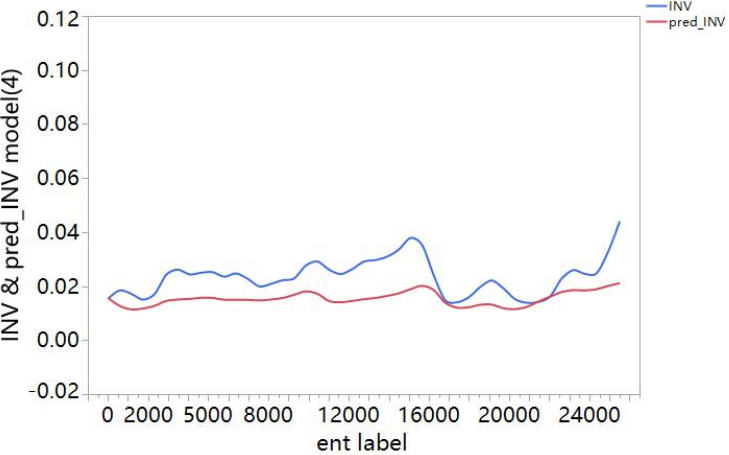}
  }
  \subfloat[Results of model(5)]
  {
      \label{fig:6}\includegraphics[width=0.5\textwidth]{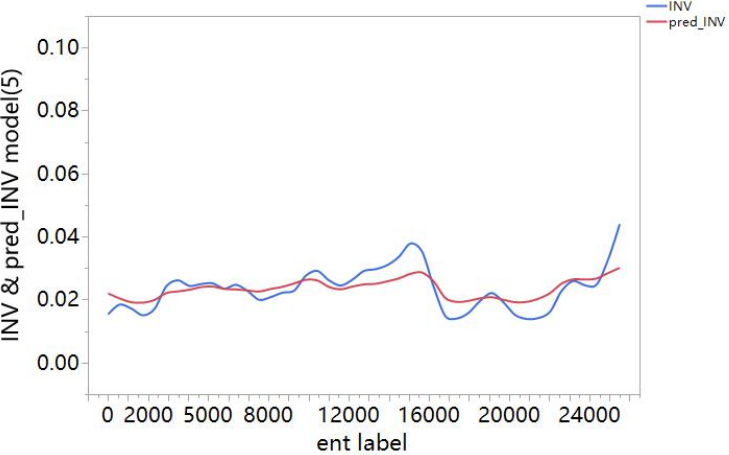}
  }
  \caption{Results of model(4) and model(5)}
  \label{training_result} 
\end{figure}

Based on the mediation effect calculation formula to analyze the mediation effect of double principal-agent costs, it is found that the mediation effect of AC1 between HOLD and INV is (-0.0422) × (-0.0247) /0.00441 = 23.6\% ; The medi ation effect of AC2 between HOLD and INV is (0.00142) × (-0.0348) /0.00441=-1.12\%, and the mediation effect of AC2 is much lower than that of AC1.


\subsection{Robustness Test}
This paper confirms that, based on the perspective of dual agency costs, HOLD can significantly reduce INV by reducing AC1, and can also significantly inhibit INV by increasing AC2. However, there may be endogenous problems in this paper, that is, it may not be because HOLD increases INV, but the benefits brought about by the increase in INV make companies tend to encourage executives with equity. In addition, there may be hysteresis in the promotion of INV brought about by executive equity incentives, and the endogenous problems in the research may lead to biased regression results test. In this paper, by replacing the index of INV with the index of the next period, after lagging one period of regression, it is found that HOLD and the company's future INV remain at the level of 1\%. HOLD can increase the company's future INV, and the double principal-agent cost still has a negative effect on the company's future INV. Confirm the robustness of the conclusion of this paper. The regression results are displayed in Table 4.

\begin{table}[H]
\caption{Robustness checks (Note: ***, ** and * imply p\textless{}0.01, p\textless{}0.05 and *p\textless{}0.1, respectively.)}\label{tab4}
\centering
\begin{tabular}{cccc}
\hline
VARIABLES & \begin{tabular}[c]{@{}c@{}}(1)\\ LINV\end{tabular} & \begin{tabular}[c]{@{}c@{}}(2)\\ LINV\end{tabular} & \begin{tabular}[c]{@{}c@{}}(3)\\ LINV\end{tabular} \\
\hline
\multirow{2}{*}{HOLD} & 0.00542*** & 0.00429*** & 0.00538*** \\
 & (7.470) & (5.932) & (7.420) \\
\multirow{2}{*}{AC1} &  & -0.0206*** &  \\
 &  & (-18.75) &  \\
\multirow{2}{*}{AC2} &  &  & -0.0311*** \\
 &  &  & (-5.478) \\
\multirow{2}{*}{LOSS} & 0.00270*** & 0.00389*** & 0.00243*** \\
 & (7.540) & (10.77) & (6.730) \\
\multirow{2}{*}{DUAL} & 0.000191 & 6.12e-05 & 0.000177 \\
 & (0.796) & (0.258) & (0.737) \\
\multirow{2}{*}{P} & 0.00525*** & 0.00488*** & 0.00519*** \\
 & (26.70) & (24.92) & (26.38) \\
\multirow{2}{*}{AGE} & -0.000186*** & -0.000196*** & -0.000179*** \\
 & (-9.680) & (-10.28) & (-9.272) \\
\multirow{2}{*}{SIZE} & -0.00108** & -0.000633*** & -0.00105*** \\
 & (-8.236) & (-4.827) & (-8.040) \\
\multirow{2}{*}{TQ} & 0.00234*** & 0.00213*** & 0.00234*** \\
 & (23.61) & (21.61) & (23.67) \\
\multirow{2}{*}{NCPS} & 0.00121*** & 0.00126*** & 0.00119*** \\
 & (9.734) & (10.29) & (9.635) \\
\multirow{2}{*}{GROWTH} & -0.00194*** & -0.00126*** & -0.00197*** \\
 & (-6.105) & (-3.967) & (-6.197) \\
\multirow{2}{*}{Constant} & -0.0344*** & -0.0427*** & -0.0336*** \\
 & (-11.34) & (-14.05) & (-11.08) \\
Year/Ind & yes & yes & yes \\
Observations & 20,581 & 20,581 & 20,581 \\
R-squared & 0.388 & 0.388 & 0.388 \\ \hline
\end{tabular}%
\end{table}

\section{Conclusion}
This paper takes China's companies listed in SH and SZ A-share market as the object to explore, and builds a model based on the principal-agent theory. From the perspective of dual agency costs, the paper analyzes the relationship between HOLD and INV. Experiment results demonstrate that HOLD can promote INV. In addition, HOLD can on one hand increase INV by reducing AC1, and on the other hand increase AC2 and thus inhibit INV. Finally, the intermediary effect of AC1 is stronger than that of AC2. The research in this paper can be deployed to guide the strategic practice of enterprise management.  

This paper verifies the mechanism of HOLD in promoting INV from the perspective of dual agency costs, but it has not carried out an in-depth analysis on the specific effects of the two types of principal-agent costs. In addition, existing studies have found that HOLD will also have an impact on external investment institutions based on the signaling theory, thereby affecting corporate financing constraints and thus affecting INV. Follow-up research can further explore the impact of HOLD on INV based on these mechanisms.

%
%
\bibliographystyle{splncs04}
\bibliography{APWEB-2023ref}

\end{document}